\documentclass{article}
\usepackage{amssymb,latexsym}
\usepackage{amsmath}
\usepackage{graphics}
\usepackage{graphicx}
\newtheorem{theorem}{Theorem}

\numberwithin{equation}{section} \numberwithin{theorem}{section}

\begin{document}
\title{The Euler-Savary Formula for One-Parameter Planar Hyperbolic Motion}
\author{Soley ERSOY, Murat TOSUN\\
         Department of Mathematics,
Faculty of Arts and Sciences
\\
Sakarya University, 54187 Sakarya/TURKEY \\}

\date{March 17, 2009}

\maketitle

\begin{abstract}
One-parameter hyperbolic planar motion was first studied by S.
Y$\ddot{\texttt{u}}$ce and N. Kuruo$\tilde{\texttt{g}}$lu.
Moreover, they analyzed the relationships between the absolute,
relative and sliding velocities of one-parameter hyperbolic planar
motion as well as the related pole curves, \cite{Yuc}.
One-parameter planar motions in the Euclidean plane $\mathbb{E}^2$
and the Euler-Savary formula in one-parameter planar motions were
given by M$\ddot{\texttt{u}}$ller, \cite{Mul}. In the present
article, one hyperbolic plane moving relative to two other
hyperbolic planes, one moving and the other fixed, was taken into
consideration and the relation between the absolute, relative and
sliding velocities of this movement was obtained. In addition, a
canonical relative system for one-parameter hyperbolic planar
motion was defined. Euler-Savary formula, which gives the
relationship between the curvature of trajectory curves, was
obtained with the help of this relative system.


\end{abstract}

\section{Preliminaries}\label{S:intro}

Before proceeding any further, we require a definition for the set
of hyperbolic number and assume the existence of any number $j$
which has the property $ j \ne  \pm 1$ . In terms of the standard
basis $\left\{ {1,j} \right\}$, the hyperbolic number  can be
written as
\[
z = x + jy
\]
where $j\left( {j^2  = 1} \right)$ is the unipotent (hyperbolic)
imaginary unit and the reel numbers $x$ and $y$ are called the
real and unipotent (or hallucinatory) parts of the hyperbolic
number $z$, respectively,
\cite{Cat}-\cite{Coc},\cite{Fje}-\cite{Fje3},\cite{Roc},
\cite{Sob}. The set of the hyperbolic numbers is
\[
 \mathbb{H}= \mathbb{R} \left[ j \right] = \left\{ {\left. {z = x + jy} \right|x,y \in ,j^2  = 1} \right\}
\]
In just the same way $ \mathbb{C}= \mathbb{R}\left[ i \right]$ are
the complex numbers extended to include the imaginary $i\quad
\left( {i^2  =  - 1} \right)$ number \cite{Yag}, the hyperbolic
numbers are the real numbers extended to include the unipotent $j$
number.\\
The hyperbolic numbers are also called perplex numbers \cite{Fje},
split-complex numbers \cite{Alf} or double numbers
\cite{Alf},\cite{Fje2},\cite{Roc}. The hyperbolic number systems
serve as the coordinates in the Lorentzian plane in the same way
as the complex numbers serve as coordinates in the Euclidean
plane. The role played by the complex numbers in Euclidean space
is played by the hyperbolic number systems in the pseudo-Euclidean
space, \cite{Sob}.\\
Addition and multiplication of the hyperbolic numbers are as
follows:
\[
\displaylines{
  \left( {x + jy} \right) + \left( {u + jv} \right) = \left( {x + u} \right) + j\left( {y + v} \right), \cr
  \left( {x + jy} \right)\left( {u + jv} \right) = \left( {xu + yv} \right) + j\left( {xv + yu} \right). \cr}
\]
This multiplication is commutative, associative and distributes
over addition. The hyperbolic conjugate of $z = x + jy$ is defined
by $\overline z  = x - jy$. The hyperbolic inner product is
\[
\left\langle {z,w} \right\rangle  = {\mathop{\rm Re}\nolimits}
\left( {z\overline w } \right) = {\mathop{\rm Re}\nolimits} \left(
{\overline z w} \right) = xu - yv
\]
where; $z = x + jy$ and $w = u + jv$. Hyperbolic numbers $z$ and
$w$ are hyperbolic (Lorentzian) orthogonal if  $ \left\langle
{z,w} \right\rangle  = 0$. Hyperbolic modulus of $z = x + jy$ is
\[
\left\| z \right\|_h  = \sqrt {\left| {\left\langle {z,z}
\right\rangle } \right|}  = \sqrt {\left| {z\overline z } \right|}
= \sqrt {\left| {x^2  - y^2 } \right|}
\]
and it is the hyperbolic distance of the point $z$ from the
origin. This is the Lorentz invariant of two-dimensional special
relativity and their unimodular multiplicative group (the group
composed of quadratic matrices determinant of which equals to 1)
is the special relativity Lorentz group, \cite{Yag2}. These
relations have been used to extend special relativity.
Furthermore, by using the functions of the hyperbolic variable,
two-dimensional special relativity has been generalized,
\cite{Cat2}. These applications make the hyperbolic numbers
appropriate for physics and the application of hyperbolic numbers
is similar to the application of complex numbers to the Euclidean
plane geometry, \cite{Yag2}. Note that the points $z \ne 0$ on the
lines $y = x$ are isotropic in the sense that they are nonzero
vectors with $ \left\| z \right\|_h  = 0$. By this way, the
hyperbolic distance creates Lorentzian geometry in $\mathbb{R}^2$.
This is different from the usual Euclidean geometry of the complex
plane, where $\left\| z \right\|_h  = 0$ only if $z=0$ in the
complex plane. The set of all points in the hyperbolic plane that
satisfy the equation $\left\| z \right\|_h  = r > 0$ is a
four-branched hyperbola of hyperbolic radius $r$, \cite{Sob}.\\
The hyperbolic number $z = x + jy$ can be written as follows:\\
While the hyperbolic number $z$ is on H-I or H-III plane, then

$$z = \pm r\left( {\cosh \varphi  + j\sinh \varphi } \right) =  \pm
re^{j\varphi },$$ While the hyperbolic number $z$ is on H-II or
H-IV plane, then
$$z = \pm r\left( {\sinh \varphi  + j\cosh \varphi
} \right) =  \pm rje^{j\varphi },$$ [See Figure1.1.]

\begin{center}
\hfil\scalebox{1}{\includegraphics{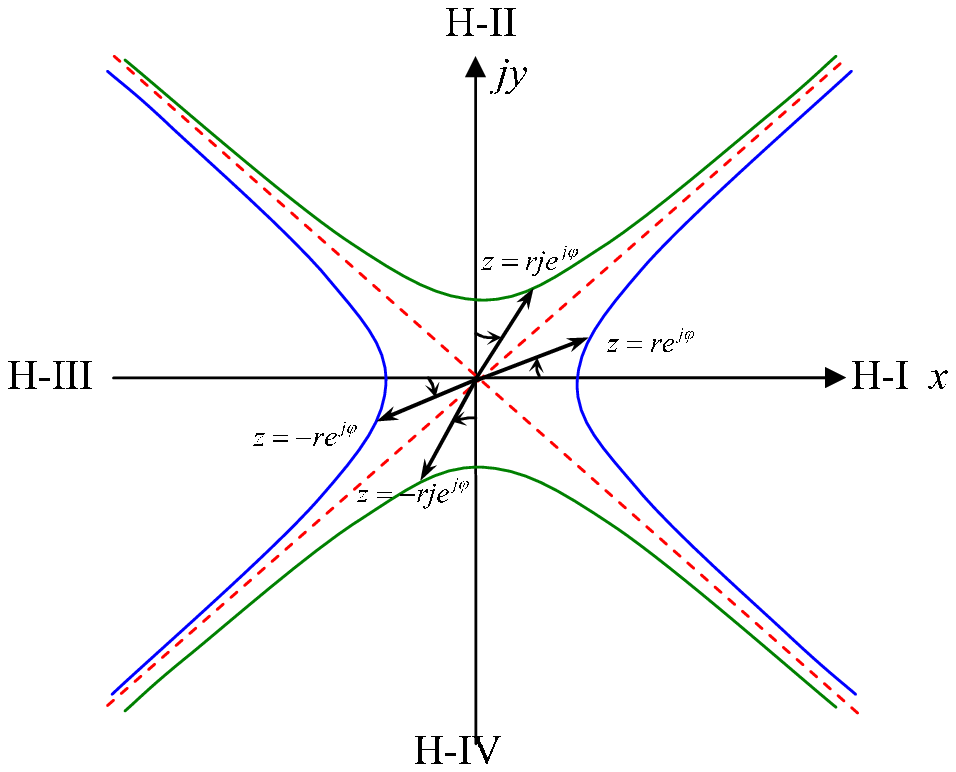}}\hfil\\
\scriptsize{Figure 1.1. Hyperbolic Plane}\\
\end{center}
\normalsize This formula can be derived by a power
series expansion due to the fact that cosh has only even powers
whereas sinh has odd powers. For all real values of the hyperbolic
angle $\varphi$, the hyperbolic number $e^{j\varphi }$ has norm 1
and lies on the right
branch of the unit hyperbola, \cite{Sob}.\\
A hyperbolic rotation defined by $e^{j\varphi }$ corresponds to
multiplication by the matrix, \cite{Sob};
\[
\left[ {\begin{array}{*{20}c}
   {\cosh \varphi } & {\sinh \varphi }  \\
   {\sinh \varphi } & {\cosh \varphi }  \\
\end{array}} \right].
\]
Another property of the hyperbolic inner product is
\[
\left\langle {ze^{j\varphi } ,we^{j\varphi } } \right\rangle  =
\left\langle {z,w} \right\rangle .
\]
In addition, a vector multiplied by $j$ is a hyperbolic orthogonal
vector, \cite{Sob}. This is similar to the role played by the
multiplication  $ i = e^{i\left( {{\pi  \mathord{\left/
 {\vphantom {\pi  2}} \right.
 \kern-\nulldelimiterspace} 2}} \right)}$ in the complex plane.

\section{Planar Hyperbolic Motion }

Let's consider an $\mathbb{A}$ plane which moves with regard to
$\mathbb{H}$ and $\mathbb{H'}$ hyperbolic planes, first one moving
and the second one fixed. Let's examine the motion of the
coordinate system  $\left\{ {B;\;{\bf a}_{\bf 1} ,\;{\bf a}_{\bf
2} } \right\}$ which defines hyperbolic plane $\mathbb{A}$, and
hyperbolic planes $\mathbb{H}$ and $\mathbb{H'}$ with regard to
the coordinate systems $\left\{ {O;\;{\bf h}_{\bf 1} ,\;{\bf
h}_{\bf 2} } \right\}$ and $\left\{ {O';\;{\bf h'}_1 ,\;{\bf h'}_2
} \right\}$. [See Figure 2.1. and 2.2.] If the vector
$\overrightarrow {OB}$ is defined by the hyperbolic number ${\bf
b} = b_1  + jb_2$, by applying the hyperbolic inner product,
$b_1^2  - b_2^2  > 0$ or $b_1^2  - b_2^2 < 0$ can be obtained. As
seen in Figure 2.1. and Figure 2.2. respectively, the vector
$\overrightarrow {OB}$ can be on the plane H-I or H-II in
hyperbolic motion.
\hfil\scalebox{1}{\includegraphics{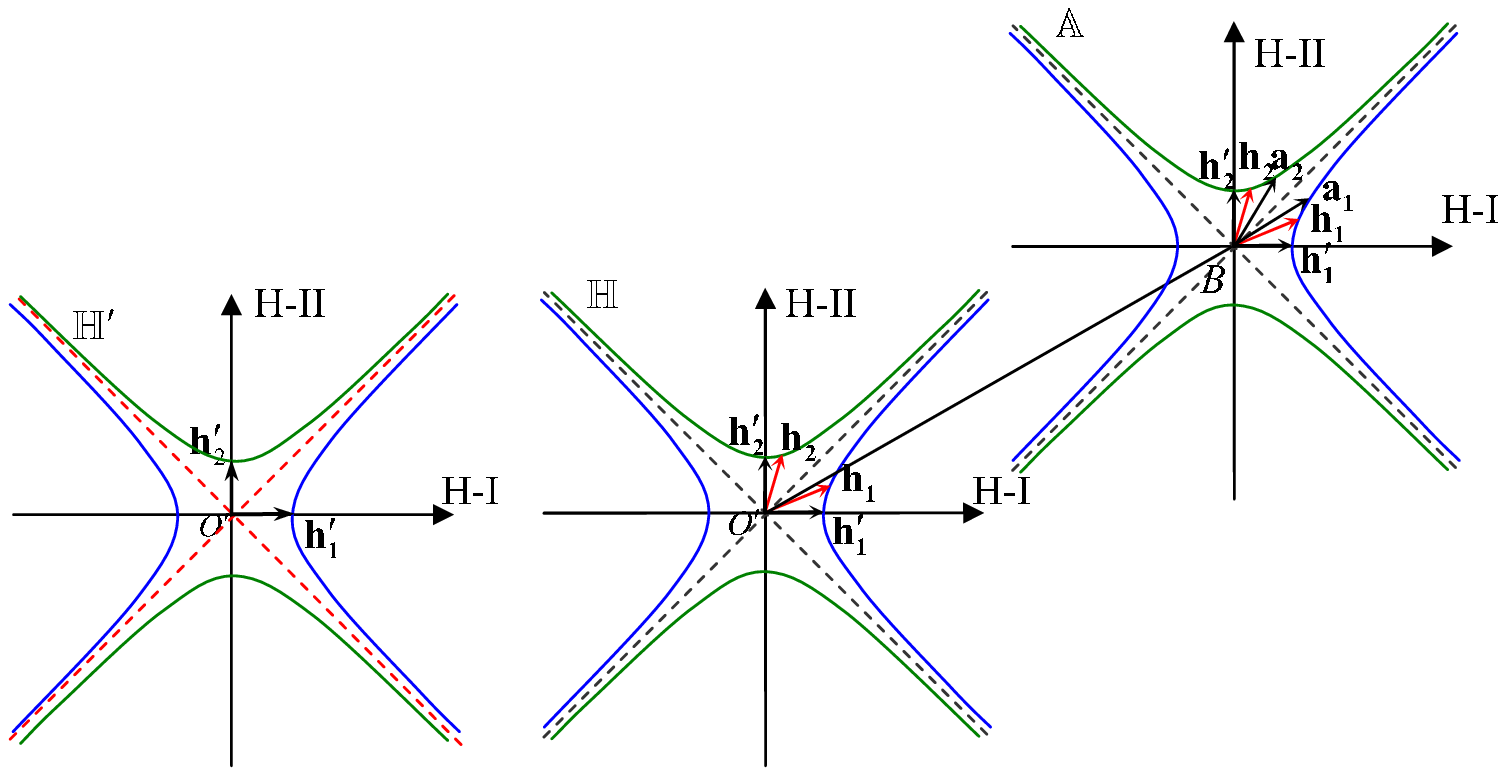}}\hfil
\begin{center}
\scriptsize{Figure2.1. $\overrightarrow {OB}$ vector is on H-I
plane}
\end{center}
\hfil\scalebox{1}{\includegraphics{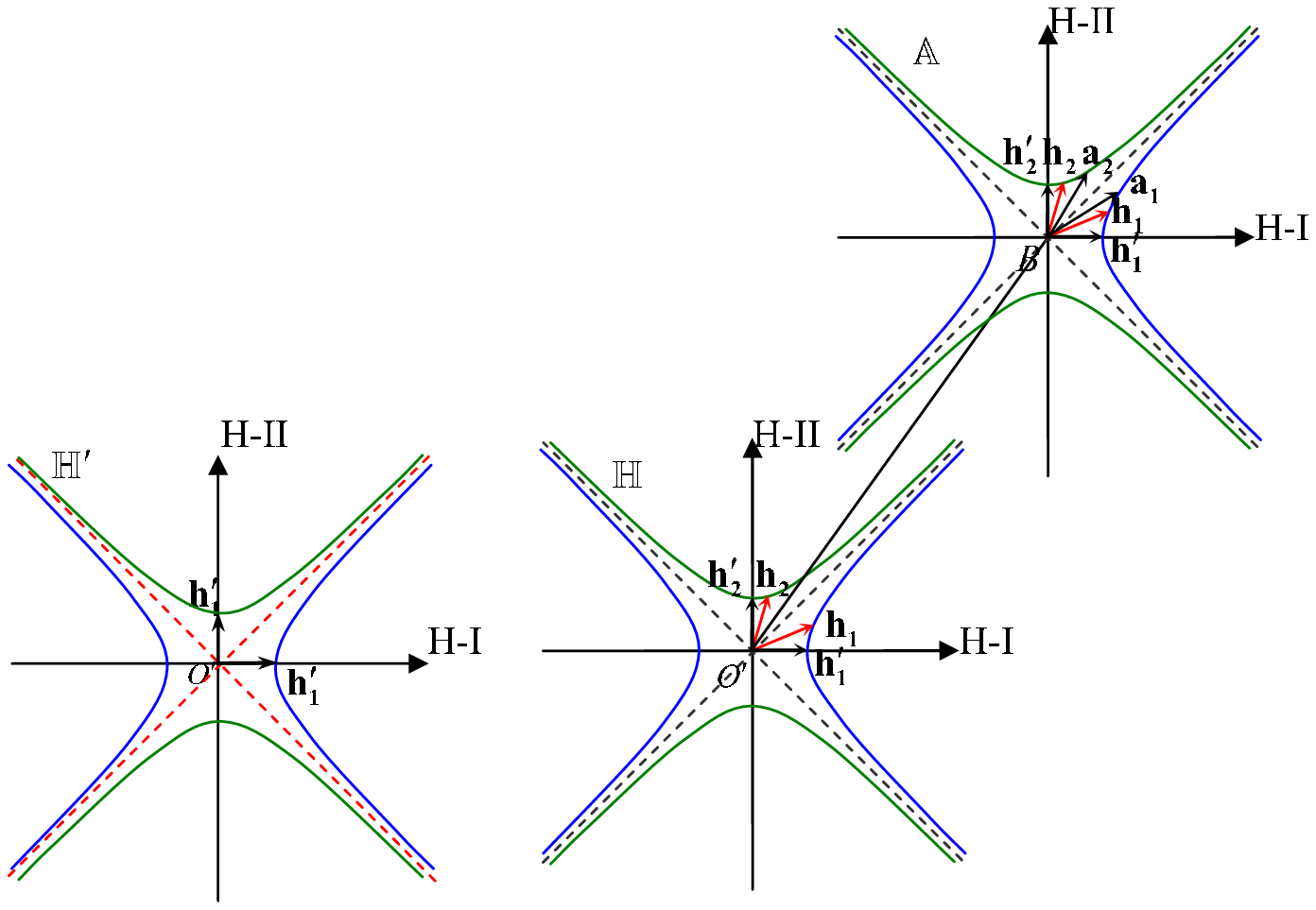}}\hfil
\begin{center}
\scriptsize{Figure2.2. $\overrightarrow {OB}$ vector is on H-II
plane}
\end{center}
\normalsize The rotation angles of the one-parameter planar
hyperbolic motion ${\mathbb{A} \mathord{\left/
 {\vphantom {A H}} \right.
 \kern-\nulldelimiterspace} \mathbb{H}}$ and  ${\mathbb{A} \mathord{\left/
 {\vphantom {A H}} \right.
 \kern-\nulldelimiterspace} \mathbb{H'}}$ are $\varphi$ and $\psi$, respectively. If the origin points of $O,\;B $ and $O',\;B $
 are coincident, then there exists following relations;
\[
\begin{array}{l}
 {\bf a}_{\bf 1}  = \cosh \varphi {\bf h}_{\bf 1}  + \sinh \varphi {\bf h}_{\bf 2}  \\
 {\bf a}_{\bf 2}  = \sinh \varphi {\bf h}_{\bf 1}  + \cosh \varphi {\bf h}_{\bf 2}  \\
 \end{array}
\]
and
\[
\begin{array}{l}
 {\bf a}_{\bf 1}  = \cosh \psi {\bf h'}_{\bf 1}  + \sinh \psi {\bf h'}_{\bf 2}  \\
 {\bf a}_{\bf 2}  = \sinh \psi {\bf h'}_{\bf 1}  + \cosh \psi {\bf h'}_{\bf 2}  \\
 \end{array}
\]
respectively [See Figure 2.3. and 2.4.].\\

\hfil\scalebox{1}{\includegraphics{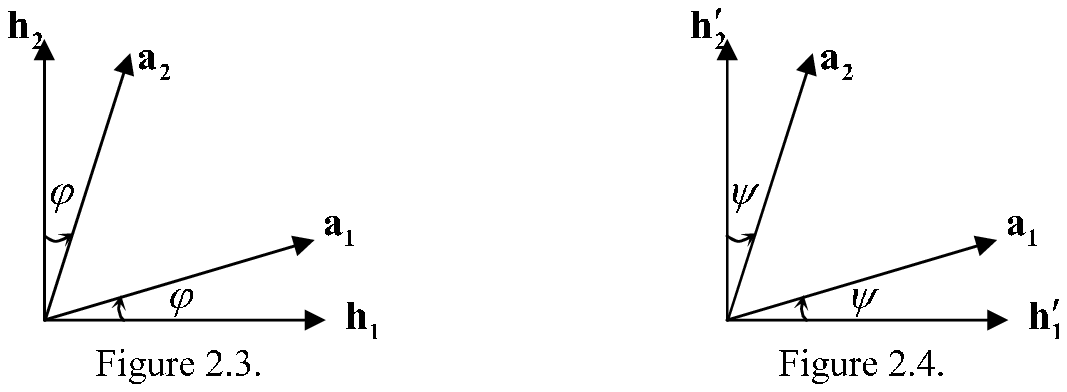}}\hfil\\
If we denote the vectors  $\overrightarrow {BX}$, $\overrightarrow
{OB}$ and  $\overrightarrow {OB'}$ with the hyperbolic numbers
${\bf{\tilde x}} = x_1  + jx_2 $, ${\bf b} = b_1  + jb_2$ and
${\bf b'}=b'_1  + jb'_2$ on the moving coordinate system of
$\mathbb{A}$, respectively; then we have
\begin{equation}\label{E:1}
{\bf{x }}= \left( {\bf{b} + \bf{\tilde x}} \right)e^{j\varphi }
\end{equation}
and
\begin{equation}\label{E:2}
\bf{x'} = \left( {\bf{b} + \bf{\tilde x}} \right)e^{j\psi }
\end{equation}
where, the hyperbolic numbers $\bf{x}$ and $\bf{x'}$ denote the
point $X$ with respect to the coordinate systems of $\mathbb{H}$
and $\mathbb{H'}$,
respectively.\\
Let's find the velocities of the one-parameter motion with the
help of the differentiation of the equations (\ref{E:1}) and
(\ref{E:2}). By differentiating equation (\ref{E:1}), we get

\begin{equation}\label{E:3}
d{\bf{x}} = \left( {\sigma  + j\tau {\bf{\tilde x}}  + d\bf{\tilde
x}} \right)e^{j\varphi }
\end{equation}
in which
\begin{equation}\label{E:4}
\sigma  = \sigma _1  + j\sigma _2  = d{\bf{b}} + j{\bf{b}}d\varphi
\quad ,\quad \tau  = d\varphi
\end{equation}
and the relative velocity vector of $X$ (with respect to
$\mathbb{H}$) is $
{\bf{V_r} } = \frac{{d{\bf{x}}}}{{dt}}$.\\
If we assume the differentiation of the equation (\ref{E:2}),
\begin{equation}\label{E:5}
d'{\bf{x}} = \left( {\sigma'  + j\tau' {\bf{\tilde x}}  +
d\bf{\tilde x}} \right)e^{j\psi}
\end{equation}
can be obtained along with the equation
\begin{equation}\label{E:6}
\sigma ' = \sigma '_1  + j\sigma '_2  = d'{\bf{b}} +
j{\bf{b'}}d\psi \quad ,\quad \tau ' = d\psi
\end{equation}
Also, the absolute velocity vector, that is, the velocity vector
of $X$ with respect to $\mathbb{H'}$, is ${\bf{V_a}}  = \frac{{d'{\bf{x}}}}{{dt}}$.\\
Here, $ \sigma _i ,\;\;\sigma '_i ,\;\;\left( {i = 1,2}
\right),\;\;\tau ,\;\;\tau ' $ are linear differential forms of
$t$ and are called Lorentzian Pfaffian forms of one-parameter
hyperbolic motion. The real parameter $t$ represents time.\\
If $\bf{V_r= 0}$ and $\bf{V_a= 0}$, the point $X$ is fixed on the
hyperbolic planes $\mathbb{H}$ and $\mathbb{H'}$, respectively.
Thus, the conditions of $X$ being fixed on the $\mathbb{H}$ and
$\mathbb{H'}$ planes are
\begin{equation}\label{E:7}
d{\bf{\tilde x}} =  - \sigma  - j\tau {\bf{\tilde x}}
\end{equation}
and
\begin{equation}\label{E:8}
d{\bf{\tilde x}} =  - \sigma ' - j\tau '{\bf{\tilde x}}
\end{equation}
respectively. If the equation (\ref{E:7}) is substituted into
equation (\ref{E:5}),
\begin{equation}\label{E:9}
d_f {\bf{x}}= \left[ {\left( {\sigma ' - \sigma } \right) +
j\left( {\tau ' - \tau } \right){\bf{\tilde x}}} \right]e^{j\psi }
\end{equation}
can be obtained, where the sliding velocity vector of the point
$X$ is  $V_f  = \frac{{d_f {\bf{x}}}}{{dt}}$. Thus, following can
be easily obtained:
\begin{equation}\label{E:10}
d'{\bf{x}} = d_f {\bf{x}} + d{\bf{x}}
\end{equation}
Just to avoid translation, it is assumed that $\dot \varphi  \ne
0$ and  $ \dot \psi  \ne 0$. The rotation pole of the motion $
{\mathbb{H} \mathord{\left/
 {\vphantom {H {H'}}} \right.
 \kern-\nulldelimiterspace} {\mathbb{H'}}}$ is characterized by the sliding velocity  $P$ being 0. For that
reason, if $d_f {\bf{x}}={\bf{0}}$, from the equation (\ref{E:9}),
the pole point $P$ of the one-parameter planar hyperbolic motion
is obtained as
\begin{equation}\label{E:11}
{\bf{p}} = j\frac{{\sigma ' - \sigma }}{{\tau  - \tau '}}
\end{equation}
and if Lorentzian coordinates are preferred on the condition that
$\overrightarrow {BP}  = {\bf{p}} = p_1  + jp_2 $, it can be
written
\begin{equation}\label{E:12}
p_1  = \frac{{\sigma '_2  - \sigma _2 }}{{\tau  - \tau '}}\quad
,\quad p_2  = \frac{{\sigma '_1  - \sigma _1 }}{{\tau  - \tau '}}
\end{equation}
which is given in \cite{Erg}.\\
In the ${\mathbb{H} \mathord{\left/
 {\vphantom {H {H'}}} \right.
 \kern-\nulldelimiterspace} {\mathbb{H'}}}$ one-parameter planar hyperbolic motion, moving and fixed
pole curves determine the geometric locus of the point $P$ in
$\mathbb{H}$ and $\mathbb{H'}$ planes, respectively. In other
words; $\left( P \right)$ and $\left( P'\right)$ are the
representation of the moving and fixed pole curves, respectively.
Also, the pole tangents can be
either on the plane H-I or H-II [See Figure 2.5. and 2.6.].\\

\hfil\scalebox{1}{\includegraphics{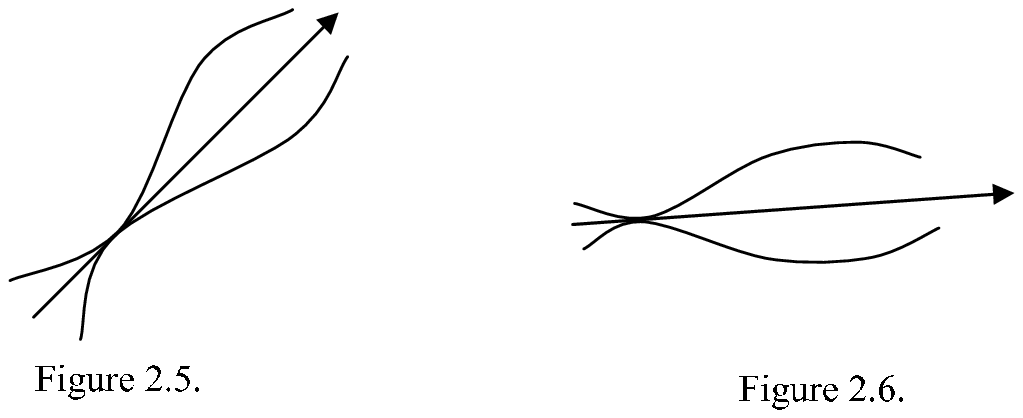}}\hfil\\
Let's first choose the pole tangents of the pole curves $\left( P
\right)$ and $\left( P'\right)$ on the plane H-II because the same
results would be obtained by following similar operations on the
plane H-I.

\section{The Euler-Savary Formula for One-Parameter Planar Hyperbolic Motion}
Let's choose the moving plane $\mathbb{A}$, represented by the
coordinate system $ \left\{ {B;{\bf a}_{\bf 1} ,{\bf a}_{\bf 2} }
\right\}$,
in such way to meet the following conditions:\\
\textbf{i)} The origin of the system $B$ coincides with the
instantaneous rotation pole $P$\\
\textbf{ii) }The axis $ \left\{ {B;{\bf a}_2 } \right\}$ is the
pole tangent, that is, it coincides with the common tangent of the
pole curves $\left( P \right)$ and $\left( P'\right)$ (on the
plane H-II) [See Figure 3.1.].\\

\hfil\scalebox{1}{\includegraphics{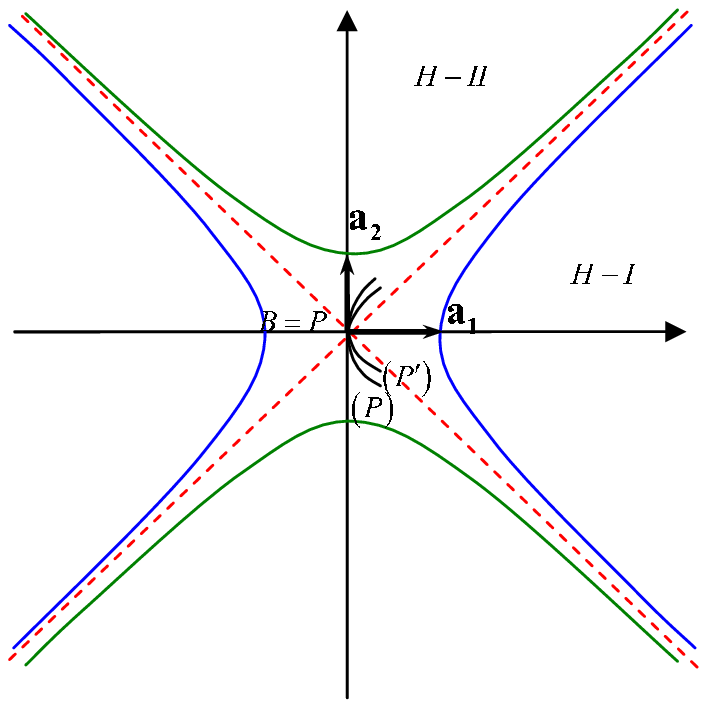}}\hfil
\begin{center}
\scriptsize{Figure 3.1.}
\end{center}
\normalsize When the condition (i) is considered: by using the
equation (\ref{E:12}),
\begin{equation}\label{E:13}
\sigma _1  = \sigma '_1 \quad ,\quad \sigma _2  = \sigma '_2
\end{equation}
are obtained. From the equations (\ref{E:4}) and  (\ref{E:6}),
\begin{equation}\label{E:14}
\begin{array}{*{20}c}
   {d{\bf{b}} = (d{\bf{b}} + j{\bf{b}}d\varphi )e^{j\varphi }  = \sigma e^{j\varphi } }  \\
   {d'{\bf{b}} = (d{\bf{b'}} + j{\bf{b'}}d\psi )e^{j\psi }  = \sigma 'e^{j\psi } }  \\
\end{array}
\end{equation}
are found. If the equation (\ref{E:13}) and the last equation are
took into consideration:
\begin{equation}\label{E:15}
d{\bf{p}} = d'{\bf{p}} = d{\bf{b}} = d'{\bf{b}}
\end{equation}
is found.\\
Thus, the moving pole curve $\left( P \right)$, the pole tangent
of which is given, and the fixed pole curves $\left( P' \right)$
are rolling on each other without sliding.\\
The second condition, that is, the condition that the pole tangent
coincides with  ${\bf a}_{\bf 2}$, requires the coefficient of
${\bf a}_{\bf 1}$ to be zero. Here, $\sigma _1  = \sigma '_1  = 0$
and $\sigma  = j\sigma _2  = j\sigma '_2$ can be written.
Consequently, the derivative equations of the canonical relative
system $\{ P;{\bf a}_{\bf 1} ,{\bf a}_{\bf 2} \}$ are
\begin{equation}\label{E:16}
d{\bf a}_1  = \tau \,{\bf a}_2  = j\tau {\kern 1pt} e^{j\varphi }
\;,\quad d{\bf a}_2  = \tau \,{\bf a}_1  =  - {\kern 1pt} \tau
{\kern 1pt} e^{j\varphi } \;,\quad d{\bf{p} }= j\sigma _2 {\kern
1pt} {\bf a}_1 = \sigma {\kern 1pt} e^{j\varphi }
\end{equation}
and
\begin{equation}\label{E:17}
d'{\bf a}_1  = \tau '{\bf a}_2  = j\tau 'e^{j\psi } ,\quad \quad
d'a_2 = \tau '{\bf a}_1  = j\tau 'e^{j\psi } ,\quad \quad
d'{\bf{p} } = j\sigma '_2 {\kern 1pt} {\bf a}_1  = \sigma {\kern
1pt} e^{j\psi }
\end{equation}
Here $\sigma  = ds$ is the scalar arc element of the pole curves
$\left( P \right)$ and $\left( P'\right)$. $\tau$ is the
hyperbolic cotangent angle, that is, two neighboring tangent
angles of $\left( P \right)$. Thus, the curvature of $\left( P
\right)$ on the point $P$ is represented by $\frac{\tau }{\sigma }
= \frac{{d\varphi }}{{ds}}$.\\
Similarly, the curvature of the hyperbolic cotangent angle $\tau'$
-that is, the fixed pole curve $\left( P' \right)$ on the point
$P$- is $\frac{{\tau '}}{\sigma } = \frac{{d\psi }}{{ds}}$.\\
The inverse values of these ratios
\begin{equation}\label{E:18}
r = \frac{\sigma }{\tau }
\end{equation}
and
\begin{equation}\label{E:19}
r' = \frac{\sigma }{{\tau '}}
\end{equation}
give the curvature radius of the pole curves $\left( P \right)$
and $\left( P'\right)$, respectively.\\
When $d\nu  = \tau ' - \tau$ is the infinitesimal small hyperbolic
instantaneous rotation angle, the moving hyperbolic plane
$\mathbb{H}$, with respect to the fixed plane $\mathbb{H'}$,
rotates around the rotation pole $P$ as much as this hyperbolic
angle in the $dt$ time scale. Thus, the hyperbolic angular
velocity of the rotational motion of $\mathbb{H}$ with respect to
$\mathbb{H'}$ is
\begin{equation}\label{E:20}
\frac{{\tau ' - \tau }}{{dt}} = \frac{{d\nu }}{{dt}} = \mathop \nu
\limits^ \bullet.
\end{equation}
From the equations (\ref{E:18}), (\ref{E:19}), and the last
equation, the following can be written:
\begin{equation}\label{E:21}
\frac{{\tau ' - \tau }}{{dt}} = \frac{{d\nu }}{{dt}} =
\frac{1}{{r'}} - \frac{1}{r}.
\end{equation}
Let the direction of the unit tangent vector ${\bf a}_{\bf 2} $ be
in the direction determined by time-based pole curves $\left( P
\right)$ and $\left( P'\right)$. Let's choose the vector ${\bf
a}_{\bf 2} $ in such way to ensure that $\frac{{ds}}{{dt}} > 0$.
In this case, $r > 0$ as the curvature center of the moving pole
$\left( P \right)$ curve is at the right side of the directed pole
tangent $\left\{ {P;{\bf a}_{\bf 2} } \right\}$. Similarly, $r' >
0$.\\
According to the canonical relative system, the differentiation
${\bf{x}}$- the coordinates of which are $x_1 ,\;x_2$- with
respect to the planes $\mathbb{H}$ and $\mathbb{H'}$ are
\begin{equation}\label{E:22}
d{\bf{x}}= \left[ {\left( {\tau x_2  + dx_1 } \right) + j(\sigma
_2  + \tau x_1  + dx_2 )} \right]e^{j\varphi }  = (\sigma  + j\tau
{\bf{x}} + d{\bf{x}})e^{j\varphi }
\end{equation}
and
\begin{equation}\label{E:23}
d'{\bf{x}} = \left[ {\left( {\tau 'x_2  + dx_1 } \right) +
j(\sigma '_2 + \tau 'x_1  + dx_2 )} \right]e^{j\varphi }  =
(\sigma  + j\tau '{\bf{x}} + d{\bf{x}})e^{j\varphi }
\end{equation}
respectively. If
\begin{equation}\label{E:24}
dx_1  = \tau x{}_2\quad {\rm and}\quad dx_2  =  - \sigma _2  -
\tau x{}_1
\end{equation}
then the point $X$ is fixed on the hyperbolic plane $\mathbb{H}$.
Similarly, if
\begin{equation}\label{E:25}
dx_1  = \tau 'x{}_2\quad {\rm and}\quad dx_2  =  - \sigma '_2  -
\tau 'x{}_1
\end{equation}
then the point $X$ is fixed on the hyperbolic plane $\mathbb{H'}$.
Also, the sliding velocity $\bf{V_f}$ of the movement ${\mathbb{H}
\mathord{\left/
 {\vphantom {H {H'}}} \right.
 \kern-\nulldelimiterspace} {\mathbb{H'}}}$ corresponds to the differentiation
\begin{equation}\label{E:26}
d_f {\bf{x}} = j\left( {\tau ' - \tau } \right)\left( {x_1  + jx_2
} \right)e^{j\varphi }  = j\left( {\tau ' - \tau }
\right){\bf{x}}e^{j\varphi }
\end{equation}
Now, let's examine the curvature centers of the trajectory curves
drawn on their fixed plane by the points of moving planes in the
motion of ${H \mathord{\left/ {\vphantom {H {H'}}} \right.
 \kern-\nulldelimiterspace} {H'}}$. In the canonical relative system, the points $X$, $X'$ having
the coordinates $x_1$, $x_2$ and $x'_1$, $x'_2$, respectively, are
situated, -together with the instantaneous rotation pole $P$ in
every $t$ moment on the instantaneous trajectory normal, which
belongs to $X$. Moreover, this curvature center can be considered
as the limit of the meeting point of the normals of the two
neighboring points on the curve. Thus,
\begin{equation}\label{E:27}
\begin{array}{l}
 \overrightarrow {PX}  = x_1  + jx_2  = {\bf{x}} \\
 \overrightarrow {PX'}  = x'_1  + jx'_2  = {\bf{x'}} \\
 \end{array}
\end{equation}
vectors have the same direction which passes through $P$. Then,
for the points $X$ and $X'$, the equation is
\begin{equation}\label{E:28}
\frac{{\bf{x}}}{{{\bf{x'}}}} = \frac{{x_1  + jx_2 }}{{x'_1  +
jx'_2 }} = \lambda \in\mathbb{R}.
\end{equation}
If the differential of this last equation is taken,
\begin{equation}\label{E:29}
\left( {x'_1 dx_1  + x'_2 dx_2  - x_1 dx'_1  - x_2 dx'_2 } \right)
+ j\left( {x'_1 dx_2  + x'_2 dx_1  - x_1 dx'_2  - x_2 dx'_1 }
\right) = 0.
\end{equation}
If the conditions that the point $X$ be fixed on the plane
$\mathbb{H}$ and the point $X'$ be fixed on the plane
$\mathbb{H'}$ are provided, then
\[
j\sigma _2 \left[ {\left( {x_1  + jx_2 } \right) - \left( {x'_1  +
jx'_2 } \right)} \right] + j\left( {x_1  + jx_2 } \right)\left(
{x'_1  + jx'_2 } \right)(\tau ' - \tau ) = 0
\]
can be obtained, that is,
\begin{equation}\label{E:30}
\sigma \left[ {{\bf{x}} - {\bf{x'}}} \right] +
j{\bf{x}}{\bf{x'}}(\tau ' - \tau ) = 0
\end{equation}
As the vectors $ \overrightarrow {PX}$, $\overrightarrow {PX'}$
are on the plane H-II ,
\begin{equation}\label{E:31}
{\bf{x}} = aje^{j\alpha }
\end{equation}
and
\begin{equation}\label{E:32}
{\bf{x'}} = a'je^{j\alpha }.
\end{equation}
That is, $a$ and $a'$, respectively, represent the distance of the
points $X$ and $X'$ on the plane H-II from the rotation pole $P$.
Also, the angle $\alpha$ is bounded by the pole curves
$\overrightarrow {PX}  = \overrightarrow {PX'}$, [See Figure 3.2.]

\hfil\scalebox{1}{\includegraphics{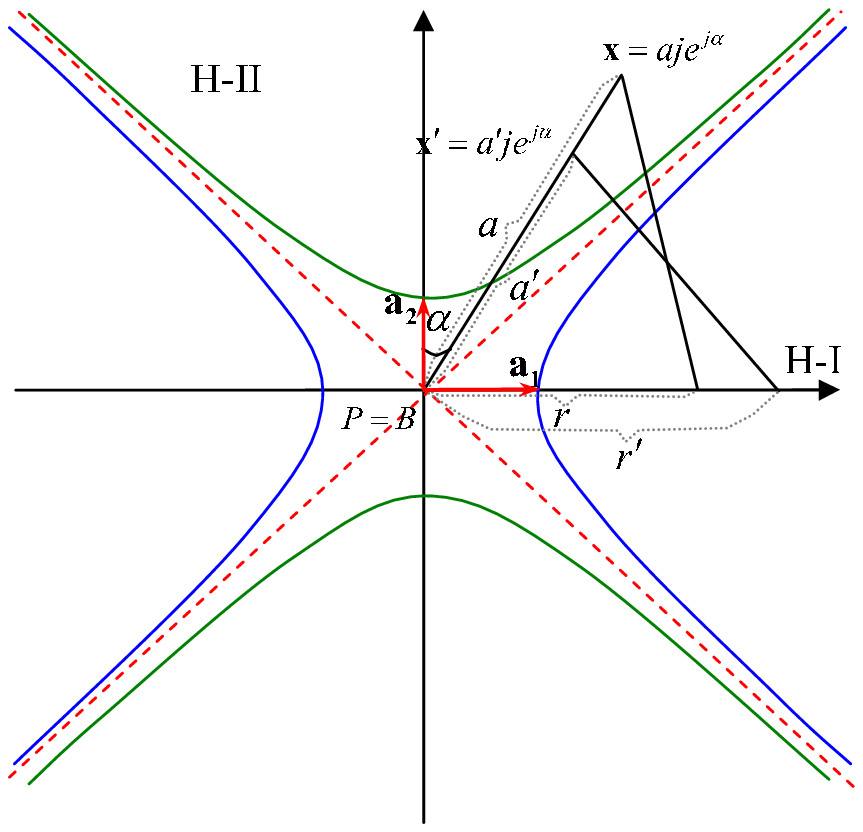}}\hfil
\begin{center}
\scriptsize{Figure 3.2.}
\end{center}
If the equations (\ref{E:31}) and (\ref{E:32}) are substituted
into equation (\ref{E:30}), then
\begin{equation}\label{E:33}
j \sigma (a - a') + jaa'e^{j\alpha } (\tau ' - \tau ) = 0
\end{equation}
can be obtained, and if the equation (\ref{E:21}) is considered
together with this last equation,
\begin{equation}\label{E:34}
\frac{{d\nu }}{{ds}} = \frac{1}{{r'}} - \frac{1}{r} = (\frac{1}{a}
- \frac{1}{{a'}})e^{ - j\alpha }
\end{equation}
is found. Here, $r$ and $r'$ are the radii of curvature of the
pole curves $P$ and $P'$, respectively. $ds$ represents the scalar
arc element and $d\nu$ represents the infinitesimal hyperbolic
angle of the motion of the pole curves.\\
The equation (\ref{E:34}) is called the Euler-Savary formula for
one-parameter plane hyperbolic motion.\\
Consequently, the following theorem can be given.
\begin{theorem}\label{T:1}
Let $\mathbb{H}$ and $\mathbb{H'}$ be the moving and fixed
hyperbolic planes, respectively. A point $X$, assumed on
$\mathbb{H}$, draws a trajectory whose instantaneous center of
curvature is $X'$ on the plane $\mathbb{H'}$ in one-parameter
planar motion ${\mathbb{H} \mathord{\left/
 {\vphantom {H {H'}}} \right.
 \kern-\nulldelimiterspace} {\mathbb{H'}}}$. In the inverse motion of ${\mathbb{H} \mathord{\left/
 {\vphantom {H {H'}}} \right.
 \kern-\nulldelimiterspace} {\mathbb{H'}}}$, a point $X'$ assumed on $\mathbb{H'}$ draws a
trajectory whose center of curvature is $X$ on the plane
$\mathbb{H}$. The relation between the points $X$ and $X'$ is
given by the Euler-Savary formula given in the equation
(\ref{E:34}).
\end{theorem}
\textbf{Remark} Let's choose the moving plane $\mathbb{A}$
represented by the coordinate system $ \left\{ {B;{\bf a}_{\bf 1}
,{\bf a}_{\bf 2} } \right\}$
in such way to meet following conditions:\\
\textbf{i)} The origin of the system $B$ and the instantaneous
rotation pole $P$ coincide with each other, i.e. $B=P$, [See
Figure 4.1.]\\
\textbf{ii)} The axis $ \left\{ {B;{\bf a}_1 } \right\}$  is the
pole tangent, that is, it coincides with the common tangent of the
pole curves $\left( P \right)$ and $\left( P'\right)$ (on the
plane H-I)

\hfil\scalebox{1}{\includegraphics{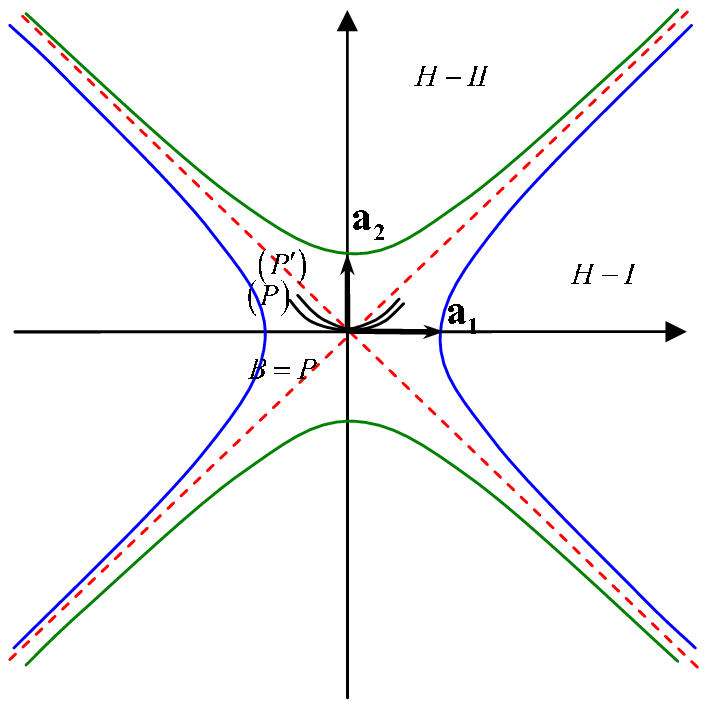}}\hfil
\begin{center}
\scriptsize{Figure 4.1.}
\end{center}
Thus, if the operations in III. section are performed considering
the conditions i) and ii), the Euler-Savary formula for
one-parameter planar hyperbolic motion remains unchanged, that is,
it is the same as in the equation (\ref{E:34})[See Figure 4.2.]

\hfil\scalebox{1}{\includegraphics{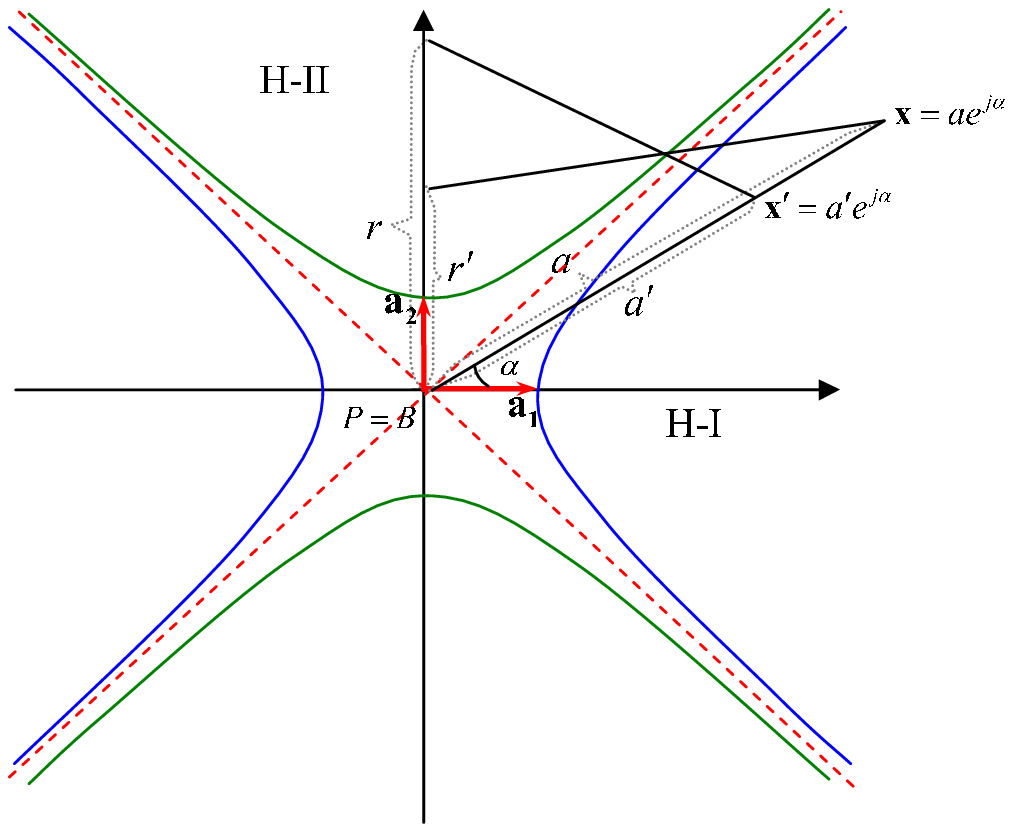}}\hfil
\begin{center}
\scriptsize{Figure 4.2.}
\end{center}

\end {document}